\documentclass{article}
\usepackage{amsmath}
\usepackage{multicol}
\usepackage{url}

\title{Calculations for OPF with Adjustable Tap Ratios}
\author{Laurens Bliek \\ Faculty of Electrical Engineering, Mathematics, and Computer Science\\ Delft University
of Technology\\ Mekelweg 4, 2628 CD Delft, The Netherlands}
\date{\today}

\newcommand{\V}{\mathcal V}

\begin{document}

\maketitle

\newpage

\tableofcontents

\newpage

\section{Introduction}

In this document, the calculations of the AC power balance and flow equations~\cite{zimmerman2010ac} used by MATPOWER will be extended to include adjustable tap ratios. The same definitions and notation as in~\cite{zimmerman2010ac} will be used, as well as the MATPOWER user manual~\cite{zimmerman2011matpower}.

\bigskip

The tap ratios are explained on pages 16 and 17 of the manual~\cite{zimmerman2011matpower}. However, in the latest version of MATPOWER (version 4.1), these tap ratios are fixed. In this document, it will be shown how to include these variables in the optimisation procedure by adding them as variables to the optimal power flow problem (OPF).

\section{Derivative of $Y_{bus}$ matrix}

In this section, the new variables $\tau$ and $\theta$ are introduced. Then the derivative of the $Y_{bus}$ matrix with respect to these new variables is given.

\bigskip

Using the definitions and the notation from~\cite{zimmerman2010ac}, we add two extra variables to the vector $X$:

$$X = \left[\begin{array}{c}\Theta\\ \V\\ P_g\\ Q_g\\ \tau\\ \theta\end{array}\right]$$

with $\tau$ and $\theta$ the transformer tap ratio magnitudes and angles respectively.

The following comes from p.\ 19 of the manual~\cite{zimmerman2011matpower}:

$$S_{bus} = [V] I_{bus}^* = [V]Y_{bus}^* V^*$$
$$Y_{bus} %= C_f^TY_f + C_t^TY_t + [Y_{sh}]
 = C_f^T[Y_{ff}]C_f + C_f^T [Y_{ft}]C_t + C_t^T[Y_{tf}]C_f + C_t^T[Y_{tt}]C_t + [Y_{sh}]$$

Note that the $Y_{bus}$ matrix is a function of the new variables $\tau$ and $\theta$, but not of any of the other variables. Using p.\ 17 of the manual~\cite{zimmerman2011matpower} for definitions of $Y_{ff}$, $Y_{ft}$, $Y_{tf}$ and $Y_{tt}$, we can calculate the following derivatives with respect to $\tau$ and $\theta$:

$$Y_{ff_\tau} = [-2Y_{ff}][\tau]^{-1}, \quad Y_{ff_\theta} = 0$$
$$Y_{ft_\tau} = [-Y_{ft}][\tau]^{-1}, \quad Y_{ft_\theta} = j[Y_{ft}]$$
$$Y_{tf_\tau} = [-Y_{tf}][\tau]^{-1}, \quad Y_{tf_\theta} = -j[Y_{tf}]$$
$$Y_{tt_\tau} = 0, \quad Y_{tt_\theta} = 0$$

Since $Y_{bus}$ is a matrix, it is more practical to calculate the derivative of $Y_{bus} \gamma$, where $\gamma$ is any vector of the right size:

\begin{align*}
(Y_{bus} \gamma)_\tau & = \frac{\partial}{\partial \tau}(C_f^T[Y_{ff}]C_f\gamma + C_f^T [Y_{ft}]C_t\gamma + C_t^T[Y_{tf}]C_f\gamma + C_t^T[Y_{tt}]C_t\gamma + [Y_{sh}]\gamma)\\
 & = \frac{\partial}{\partial \tau}(C_f^T[C_f\gamma]Y_{ff} + C_f^T [C_t\gamma]Y_{ft} + C_t^T[C_f\gamma]Y_{tf} + C_t^T[C_t\gamma]Y_{tt} + [Y_{sh}]\gamma)\\
 & = C_f^T[C_f\gamma][-2Y_{ff}][\tau]^{-1} + C_f^T[C_t\gamma][-Y_{ft}][\tau]^{-1}+C_t^T[C_f\gamma][-Y_{tf}][\tau]^{-1}
\end{align*}
\begin{align*}
(Y_{bus} \gamma)_\theta & = \frac{\partial}{\partial \theta}(C_f^T[C_f\gamma]Y_{ff} + C_f^T [C_t\gamma]Y_{ft} + C_t^T[C_f\gamma]Y_{tf} + C_t^T[C_t\gamma]Y_{tt} + [Y_{sh}]\gamma)\\
 & = C_f^T[C_t\gamma][jY_{ft}]+C_t^T[C_f\gamma][-jY_{tf}]
\end{align*}

\section{First and second derivatives of power balance equations}

The function in the power balance equation $G^s(X) = S_{bus} + S_d - C_g S_g = 0$, with $S_{bus} = [V] I_{bus}^* = [V] Y_{bus}^* V^*$, now gets two new first derivatives:

$$G_X^s = \left[G_\Theta^s\ G_\V^s\ G_{P_g}^s\ G_{Q_g}^s\ G_\tau^s\ G_\theta^s\right]$$

Using the derivatives of $Y_{bus}$ above, we have the following result for the first derivatives:

\begin{align*}
G_\tau^s  & = [V](Y_{bus}^* V^*)_\tau = [V]\Big(C_f^T[C_fV^*][-2Y_{ff}^*][\tau]^{-1} + C_f^T[C_tV^*][-Y_{ft}^*][\tau]^{-1}+C_t^T[C_fV^*][-Y_{tf}^*][\tau]^{-1}\Big)\\
G_\theta^s & = [V](Y_{bus}^* V^*)_\theta = [V]\Big(C_f^T[C_tV^*][-jY_{ft}^*]+C_t^T[C_fV^*][jY_{tf}^*]\Big)
\end{align*}

The second derivative matrix becomes a matrix consisting of $16$ terms:

$$G_{XX}^s(\lambda) = \left[\begin{array}{cccccc}
  	G_{\Theta\Theta}^s(\lambda) & G_{\Theta\V}^s(\lambda) & 0 & 0 & G_{\Theta\tau}^s(\lambda) & G_{\Theta\theta}^s(\lambda)\\
  	G_{\V\Theta}^s(\lambda) & G_{\V\V}^s(\lambda) & 0 & 0 & G_{\V\tau}^s(\lambda) & G_{\V\theta}^s(\lambda)\\
  	0 & 0 & 0 & 0 & 0 & 0\\
  	0 & 0 & 0 & 0 & 0 & 0\\
  	G_{\tau\Theta}^s(\lambda) & G_{\tau\V}^s(\lambda) & 0 & 0 & G_{\tau\tau}^s(\lambda) & G_{\tau\theta}^s(\lambda)\\
  	G_{\theta\Theta}^s(\lambda) & G_{\theta\V}^s(\lambda) & 0 & 0 & G_{\theta\tau}^s(\lambda) & G_{\theta\theta}^s(\lambda)\\
  	\end{array}\right]$$

The first four second derivatives have already been derived~\cite{zimmerman2010ac}. The derivation of the twelve new second derivatives is shown on the following pages.

\begin{align*}
G_{\Theta\tau}^s(\lambda) & = \frac{\partial}{\partial \tau}(G_\Theta^{s^T} \lambda)\\
& = \frac{\partial}{\partial \tau} j([I_{bus}^*] - [V^*]Y_{bus}^{*^T})[V]\lambda\\
& = \frac{\partial}{\partial \tau} j([V][\lambda]I_{bus}^* - [V^*]Y_{bus}^{*^T}[V]\lambda)\\
& = \frac{\partial}{\partial \tau} j([V][\lambda](C_f^T[Y_{ff}^*]C_fV^* + C_f^T [Y_{ft}^*]C_tV^* + C_t^T[Y_{tf}^*]C_fV^* + C_t^T[Y_{tt}^*]C_tV^*)\\
& \quad - [V^*](C_f^T[Y_{ff}^*]C_f + C_t^T[Y_{ft}^*]C_f^T + C_f^T[Y_{tf}^*]C_t + C_t^T[Y_{tt}^*]C_t)[V]\lambda))\\
& = \frac{\partial}{\partial \tau} j\Big([V][\lambda](C_f^T[C_fV^*]Y_{ff}^* + C_f^T [C_tV^*]Y_{ft}^* + C_t^T[C_fV^*]Y_{tf}^* + C_t^T[C_tV^*]Y_{tt}^*)\\
& \quad - [V^*](C_f^T[C_f[V]\lambda]Y_{ff}^* + C_t^T[C_f[V]\lambda]Y_{ft}^* + C_f^T[C_t[V]\lambda]Y_{tf}^* + C_t^T[C_t[V]\lambda]Y_{tt}^*)\Big)\\
& = j\Big([V][\lambda]\big(C_f^T[C_fV^*][-2Y_{ff}^*] + C_f^T [C_tV^*][-Y_{ft}^*] + C_t^T[C_fV^*][-Y_{tf}^*]\big)\\
& \quad - [V^*]\big(C_f^T[C_f[V]\lambda][-2Y_{ff}^*] + C_t^T[C_f[V]\lambda][-Y_{ft}^*] + C_f^T[C_t[V]\lambda][-Y_{tf}^*])\Big)[\tau]^{-1}\\
\\
G_{\Theta\theta}^s(\lambda) & = \frac{\partial}{\partial \theta} j\Big([V][\lambda](C_f^T[C_fV^*]Y_{ff}^* + C_f^T [C_tV^*]Y_{ft}^* + C_t^T[C_fV^*]Y_{tf}^* + C_t^T[C_tV^*]Y_{tt}^*)\\
& \quad - [V^*](C_f^T[C_f[V]\lambda]Y_{ff}^* + C_t^T[C_f[V]\lambda]Y_{ft}^* + C_f^T[C_t[V]\lambda]Y_{tf}^* + C_t^T[C_t[V]\lambda]Y_{tt}^*)\Big)\\
& = j\Big([V][\lambda](C_f^T [C_tV^*][-jY_{ft}^*] + C_t^T[C_fV^*][jY_{tf}^*])\\
& \quad - [V^*](C_t^T[C_f[V]\lambda][-jY_{ft}^*] + C_f^T[C_t[V]\lambda][jY_{tf}^*])\Big)\\
\\
G_{\V\tau}^s(\lambda) & = \frac{\partial}{\partial \tau}(G_\V^{s^T} \lambda)\\
& = \frac{\partial}{\partial \tau} [\V]^{-1}([I_{bus}^*] + [V^*]Y_{bus}^{*^T})[V]\lambda\\
& = \frac{\partial}{\partial \tau} [\V]^{-1}([V][\lambda]I_{bus}^* + [V^*]Y_{bus}^{*^T}[V]\lambda)\\
& = \frac{\partial}{\partial \tau} [\V]^{-1}([V][\lambda](C_f^T[Y_{ff}^*]C_fV^* + C_f^T [Y_{ft}^*]C_tV^* + C_t^T[Y_{tf}^*]C_fV^* + C_t^T[Y_{tt}^*]C_tV^*)\\
& \quad + [V^*](C_f^T[Y_{ff}^*]C_f + C_t^T[Y_{ft}^*]C_f^T + C_f^T[Y_{tf}^*]C_t + C_t^T[Y_{tt}^*]C_t)[V]\lambda)\\
& = \frac{\partial}{\partial \tau} [\V]^{-1}\Big([V][\lambda](C_f^T[C_fV^*]Y_{ff}^* + C_f^T [C_tV^*]Y_{ft}^* + C_t^T[C_fV^*]Y_{tf}^* + C_t^T[C_tV^*]Y_{tt}^*)\\
& \quad + [V^*](C_f^T[C_f[V]\lambda]Y_{ff}^* + C_t^T[C_f[V]\lambda]Y_{ft}^* + C_f^T[C_t[V]\lambda]Y_{tf}^* + C_t^T[C_t[V]\lambda]Y_{tt}^*)\Big)\\
& = [\V]^{-1}\Big([V][\lambda](C_f^T[C_fV^*][-2Y_{ff}^*][\tau]^{-1} + C_f^T [C_tV^*][-Y_{ft}^*][\tau]^{-1} + C_t^T[C_fV^*][-Y_{tf}^*][\tau]^{-1})\\
& \quad + [V^*](C_f^T[C_f[V]\lambda][-2Y_{ff}^*][\tau]^{-1} + C_t^T[C_f[V]\lambda][-Y_{ft}^*][\tau]^{-1} + C_f^T[C_t[V]\lambda][-Y_{tf}^*][\tau]^{-1})\Big)\\
\\
G_{\V\theta}^s(\lambda) & = \frac{\partial}{\partial \theta} [\V]^{-1}\Big([V][\lambda](C_f^T[C_fV^*]Y_{ff}^* + C_f^T [C_tV^*]Y_{ft}^* + C_t^T[C_fV^*]Y_{tf}^* + C_t^T[C_tV^*]Y_{tt}^*)\\
& \quad + [V^*](C_f^T[C_f[V]\lambda]Y_{ff}^* + C_t^T[C_f[V]\lambda]Y_{ft}^* + C_f^T[C_t[V]\lambda]Y_{tf}^* + C_t^T[C_t[V]\lambda]Y_{tt}^*)\Big)\\
& = [\V]^{-1}\Big([V][\lambda](C_f^T [C_tV^*][-jY_{ft}^*] + C_t^T[C_fV^*][jY_{tf}^*])\\
& \quad + [V^*](C_t^T[C_f[V]\lambda][-jY_{ft}^*] + C_f^T[C_t[V]\lambda][jY_{tf}^*])\Big)\\
\end{align*}

\begin{align*}
G_{\tau\Theta}^s(\lambda) &= \frac{\partial}{\partial \Theta}(G_\tau^{s^T} \lambda)\\
& = \frac{\partial}{\partial \Theta}[\tau]^{-1}\Big([-2Y_{ff}^*][C_fV^*]C_f + [-Y_{ft}^*][C_tV^*]C_f + [-Y_{tf}^*][C_fV^*]C_t\Big)[V]\lambda\\
& = \frac{\partial}{\partial \Theta}[\tau]^{-1}\Big([-2Y_{ff}^*][C_fV^*]C_f + [-Y_{ft}^*][C_tV^*]C_f + [-Y_{tf}^*][C_fV^*]C_t\Big)[\lambda]V\\
& = [\tau]^{-1}\Big([-2Y_{ff}^*][C_fV^*]C_f + [-Y_{ft}^*][C_tV^*]C_f + [-Y_{tf}^*][C_fV^*]C_t\Big)[\lambda][jV]\\
& \quad + [\tau]^{-1}\Big([-2Y_{ff}^*][C_f[\lambda]V]C_f[-jV^*] + [-Y_{ft}^*][C_f[\lambda]V]C_t[-jV^*] + [-Y_{tf}^*][C_t[\lambda]V]C_f[-jV^*]\Big)\\
& = G_{\Theta\tau}^s(\lambda)^T\\
\\
G_{\tau\V}^s(\lambda) & = \frac{\partial}{\partial \V}[\tau]^{-1}\Big([-2Y_{ff}^*][C_fV^*]C_f + [-Y_{ft}^*][C_tV^*]C_f + [-Y_{tf}^*][C_fV^*]C_t\Big)[\lambda]V\\
& = [\tau]^{-1}\Big([-2Y_{ff}^*][C_fV^*]C_f + [-Y_{ft}^*][C_tV^*]C_f + [-Y_{tf}^*][C_fV^*]C_t\Big)[\lambda][V][\V]^{-1}\\
& \quad + [\tau]^{-1}\Big([-2Y_{ff}^*][C_f[\lambda]V]C_f[V^*][\V]^{-1} + [-Y_{ft}^*][C_f[\lambda]V]C_t[V^*][\V]^{-1} \\
& \quad + [-Y_{tf}^*][C_t[\lambda]V]C_f[V^*][\V]^{-1}\Big)\\
& = G_{\V\tau}^s(\lambda)^T\\
\\
G_{\theta\Theta}^s(\lambda) & = \frac{\partial}{\partial \Theta}(G_\theta^{s^T} \lambda)\\
& = \frac{\partial}{\partial \Theta} ([-jY_{ft}^*][C_tV^*]C_f + [jY_{tf}^*][C_fV^*]C_t)[\lambda] V\\
& = \Big([-jY_{ft}^*][C_tV^*]C_f + [jY_{tf}^*][C_fV^*]C_t\Big)[\lambda][jV]\\
& \quad + [-jY_{ft}^*][C_f[\lambda] V]C_t[-jV^*] + [jY_{tf}^*][C_t[\lambda] V]C_f[-jV^*]\\
& = G_{\Theta\theta}^s(\lambda)^T\\
\\
G_{\theta\V}^s(\lambda) & = \frac{\partial}{\partial \V} ([-jY_{ft}^*][C_tV^*]C_f + [jY_{tf}^*][C_fV^*]C_t)[\lambda] V\\
& = \Big([-jY_{ft}^*][C_tV^*]C_f + [jY_{tf}^*][C_fV^*]C_t\Big)[\lambda][V][\V]^{-1}\\
& \quad + [-jY_{ft}^*][C_f[\lambda] V]C_t[V^*][\V]^{-1} + [jY_{tf}^*][C_t[\lambda] V]C_f[V^*][\V]^{-1}\\
& = G_{\V\theta}^s(\lambda)^T
\end{align*}

\begin{align*}
G_{\tau\tau}^s(\lambda) &= \frac{\partial}{\partial \tau}(G_\tau^{s^T} \lambda)\\
& = \frac{\partial}{\partial \tau} [\tau]^{-1}\Big([-2Y_{ff}^*][C_fV^*]C_f + [-Y_{ft}^*][C_tV^*]C_f + [-Y_{tf}^*][C_fV^*]C_t\Big)[V]\lambda\\
& = \frac{\partial}{\partial \tau} \Big([-2Y_{ff}^*][C_fV^*][C_f[V]\lambda] + [-Y_{ft}^*][C_tV^*][C_f[V]\lambda] + [-Y_{tf}^*][C_fV^*][C_t[V]\lambda]\Big) \tau.^{-1}\\
& = [\tau]^{-1}\Big([C_fV^*][C_f[V]\lambda][4Y_{ff}^*][\tau]^{-1} + [C_tV^*][C_f[V]\lambda][Y_{ft}^*][\tau]^{-1} + [C_fV^*][C_t[V]\lambda][Y_{tf}^*][\tau]^{-1}\Big)\\
& \quad + \Big([2Y_{ff}^*][C_fV^*][C_f[V]\lambda] + [Y_{ft}^*][C_tV^*][C_f[V]\lambda] + [Y_{tf}^*][C_fV^*][C_t[V]\lambda]\Big)[\tau]^{-2}\\
& = [\tau]^{-2}\Big([6Y_{ff}^*][C_fV^*][C_f[V]\lambda] + [2Y_{ft}^*][C_tV^*][C_f[V]\lambda] + [2Y_{tf}^*][C_fV^*][C_t[V]\lambda]\Big)\\
%\end{align*}
%\begin{align*}
\\
G_{\tau\theta}^s(\lambda) &= \frac{\partial}{\partial \theta} ([-2Y_{ff}^*][C_fV^*][C_f[V]\lambda] + [-Y_{ft}^*][C_tV^*][C_f[V]\lambda] + [-Y_{tf}^*][C_fV^*][C_t[V]\lambda]) \tau.^{-1}\\
&= \frac{\partial}{\partial \theta} [\tau]^{-1}(-2[C_fV^*][C_f[V]\lambda]Y_{ff}^* + -[C_tV^*][C_f[V]\lambda]Y_{ft}^* + -[C_fV^*][C_t[V]\lambda]Y_{tf}^*)\\
& = [\tau]^{-1}([C_tV^*][C_f[V]\lambda][jY_{ft}^*] + [C_fV^*][C_t[V]\lambda][-jY_{tf}^*])\\
%\end{align*}
%\begin{align*}
\\
G_{\theta\tau}^s(\lambda) & = \frac{\partial}{\partial \tau}(G_\theta^{s^T} \lambda)\\
& = \frac{\partial}{\partial \tau} ([-jY_{ft}^*][C_tV^*]C_f + [jY_{tf}^*][C_fV^*]C_t)[V]\lambda\\
& = \frac{\partial}{\partial \tau} (-j[C_tV^*][C_f[V]\lambda]Y_{ft}^* + j[C_fV^*][C_t[V]\lambda]Y_{tf}^*)\\
& = -j[C_tV^*][C_f[V]\lambda][-Y_{ft}^*][\tau]^{-1} + j[C_fV^*][C_t[V]\lambda][-Y_{tf}^*][\tau]^{-1}\\
& = G_{\tau\theta}^s(\lambda) = G_{\tau\theta}^s(\lambda)^T\\
\\
G_{\theta\theta}^s(\lambda) & = \frac{\partial}{\partial \theta} (-j[C_tV^*][C_f[V]\lambda]Y_{ft}^* +j[C_fV^*][C_t[V]\lambda]Y_{tf}^*)\\
& = -[C_tV^*][C_f[V]\lambda][Y_{ft}^*] -[C_fV^*][C_t[V]\lambda][Y_{tf}^*]
\end{align*}

Here, the notation $\tau.^{-1}$ is used for elementwise exponentiation of a vector.

\newpage
\section{First and second derivatives of line flow constraints}

The functions of the line flow constraints $h_f(\Theta,\V) = |F_f(\Theta,\V)| - F_{max} \leq 0$ and $h_t(\Theta,\V) = |F_t(\Theta,\V)| - F_{max} \leq 0$ get new derivatives depending on the chosen form of flow constraints (p.\ 30 of~\cite{zimmerman2011matpower}). In this document, it is assumed that the constraints depend on the current $I$, not on the real or apparent power. This gives 
$F_f = I_f$ and $F_t = I_t$, and we have the following derivatives with respect to variables $\tau$ and $\theta$:

\begin{equation}
	\frac{\partial I_f}{\partial X} = \left[\frac{\partial I_f}{\partial \Theta} \ \frac{\partial I_f}{\partial \V} \ \frac{\partial I_f}{\partial P_g} \ \frac{\partial I_f}{\partial Q_g} \ \frac{\partial I_f}{\partial \tau} \ \frac{\partial I_f}{\partial \theta}\right].
\end{equation}

The first four derivatives have already been calculated~\cite{zimmerman2010ac}. The other two can be derived as follows, for both $I_f$ and $I_t$:

\begin{align*}
	\frac{\partial I_f}{\partial \tau} & = \frac{\partial}{\partial \tau} ([Y_{ff}]C_f V + [Y_{ft}] C_t V)\\
	& = \frac{\partial}{\partial \tau} ( [C_f V]Y_{ff} + [C_t V]Y_{ft} )\\
	& = [C_f V][-2Y_{ff}][\tau]^{-1} + [C_t V][-Y_{ft}][\tau]^{-1}\\
	\frac{\partial I_f}{\partial \theta} & = \frac{\partial}{\partial \theta} ( [C_f V]Y_{ff} + [C_t V]Y_{ft} )\\
	& = [C_t V][jY_{ft}].\\
%\end{align*}
%
%Similarly,
%
%\begin{align*}
	\frac{\partial I_t}{\partial \tau} & = \frac{\partial}{\partial \tau} ([Y_{tf}]C_f V + [Y_{tt}] C_t V)\\
	& = \frac{\partial}{\partial \tau} ( [C_f V]Y_{tf} + [C_t V]Y_{tt} )\\
	& = [C_f V][-Y_{tf}][\tau]^{-1}\\
	\frac{\partial I_t}{\partial \theta} & = \frac{\partial}{\partial \theta} ( [C_f V]Y_{tf} + [C_t V]Y_{tt} )\\
	& = [C_f V][-jY_{tf}].
\end{align*}

The second derivatives matrix becomes a matrix consisting of $16$ terms, of which the first four have already been derived~\cite{zimmerman2010ac}:

$$I_{f_{XX}}(\mu) = \left[\begin{array}{cccccc}
  	I_{f_{\Theta\Theta}}(\mu) & I_{f_{\Theta\V}}(\mu) & 0 & 0 & I_{f_{\Theta\tau}}(\mu) & I_{f_{\Theta\theta}}(\mu)\\
  	I_{f_{\V\Theta}}(\mu) & I_{f_{\V\V}}(\mu) & 0 & 0 & I_{f_{\V\tau}}(\mu) & I_{f_{\V\theta}}(\mu)\\
  	0 & 0 & 0 & 0 & 0 & 0\\
  	0 & 0 & 0 & 0 & 0 & 0\\
  	I_{f_{\tau\Theta}}(\mu) & I_{f_{\tau\V}}(\mu) & 0 & 0 & I_{f_{\tau\tau}}(\mu) & I_{f_{\tau\theta}}(\mu)\\
  	I_{f_{\theta\Theta}}(\mu) & I_{f_{\theta\V}}(\mu) & 0 & 0 & I_{f_{\theta\tau}}(\mu) & I_{f_{\theta\theta}}(\mu)\\
  	\end{array}\right].$$
  	
The other $12$ derivatives can be derived as follows:

\begin{align*}
I_{f_{\Theta\tau}}(\mu) & = \frac{\partial}{\partial \tau}(I_{f_\Theta}^T \mu)\\
	& = \frac{\partial}{\partial \tau}(j [V](C_f^T [Y_{ff}] + C_t^T [Y_{ft}]) \mu)\\
	& = \frac{\partial}{\partial \tau}(j [V](C_f^T [\mu] Y_{ff} + C_t^T [\mu]Y_{ft}))\\
	& = j [V](C_f^T [\mu] [-2Y_{ff}][\tau]^{-1} + C_t^T [\mu][-Y_{ft}][\tau]^{-1})\\
I_{f_{\Theta\theta}}(\mu) & = \frac{\partial}{\partial \theta}(j [V](C_f^T [\mu] Y_{ff} + C_t^T [\mu]Y_{ft}))\\
	& = j [V]C_t^T [\mu][jY_{ft}]\\
I_{f_{\V\tau}}(\mu) & = \frac{\partial}{\partial \tau}(I_{f_\V}^T \mu)\\
	& = \frac{\partial}{\partial \tau}\Big([\V]^{-1}[V](C_f^T [Y_{ff}] + C_t^T [Y_{ft}]) \mu\Big)\\
	& = \frac{\partial}{\partial \tau}\Big([\V]^{-1}[V](C_f^T [\mu] Y_{ff} + C_t^T [\mu]Y_{ft})\Big)\\
	& = [\V]^{-1}[V](C_f^T [\mu] [-2Y_{ff}][\tau]^{-1} + C_t^T [\mu][-Y_{ft}][\tau]^{-1})\\
I_{f_{\V\theta}}(\mu) & = \frac{\partial}{\partial \theta}\Big([\V]^{-1}[V](C_f^T [\mu] Y_{ff} + C_t^T [\mu]Y_{ft})\Big)\\
	& = [\V]^{-1}[V]C_t^T [\mu][jY_{ft}]\\
\\
I_{f_{\tau\Theta}}(\mu) & = \frac{\partial}{\partial \Theta}(I_{f_\tau}^T \mu)\\
	& = \frac{\partial}{\partial \Theta}([\tau]^{-1}[-2Y_{ff}][C_fV]\mu + [\tau]^{-1}[-Y_{ft}][C_tV]\mu)\\
	& = [\tau]^{-1}[-2Y_{ff}][\mu]C_f[jV] + [\tau]^{-1}[-Y_{ft}][\mu]C_t[jV]\\
	& = I_{f_{\Theta\tau}}(\mu)^T\\
I_{f_{\tau\V}}(\mu) & = \frac{\partial}{\partial \V}([\tau]^{-1}[-2Y_{ff}][C_fV]\mu + [\tau]^{-1}[-Y_{ft}][C_tV]\mu)\\
 & = [\tau]^{-1}[-2Y_{ff}][\mu]C_f[V][\V]^{-1} + [\tau]^{-1}[-Y_{ft}][\mu]C_t[V][\V]^{-1}\\
 & = I_{f_{\V\tau}}(\mu)^T\\
I_{f_{\theta\Theta}}(\mu) & = \frac{\partial}{\partial \Theta}(I_{f_\theta}^T \mu)\\
	& = \frac{\partial}{\partial \Theta}([jY_{ft}][C_t V]\mu)\\
	& = [jY_{ft}][\mu]C_t [jV]\\
	& = I_{f_{\Theta\theta}}(\mu)^T\\
I_{f_{\theta\V}}(\mu) & = \frac{\partial}{\partial \V}([jY_{ft}][C_t V]\mu)\\
	& = [jY_{ft}][\mu]C_t [V][V]^{-1}\\
	& = I_{f_{\V\theta}}(\mu)^T\\
\end{align*}
\begin{align*}
I_{f_{\tau\tau}}(\mu) & = \frac{\partial}{\partial \tau}([\tau]^{-1}[-2Y_{ff}][C_fV]\mu + [\tau]^{-1}[-Y_{ft}][C_tV]\mu)\\
	& = \frac{\partial}{\partial \tau}(-2[\tau]^{-1}[C_fV][\mu]Y_{ff} - [\tau]^{-1}[C_tV][\mu]Y_{ft})\\
	& = 6[\tau]^{-2}[C_fV][\mu][Y_{ff}] + 2[\tau]^{-2}[C_tV][\mu][Y_{ft}])\\
I_{f_{\tau\theta}}(\mu) & = \frac{\partial}{\partial \theta}(-2[\tau]^{-1}[C_fV][\mu]Y_{ff} - [\tau]^{-1}[C_tV][\mu]Y_{ft})\\
	& = - [\tau]^{-1}[C_tV][\mu][jY_{ft}]\\
I_{f_{\theta\tau}}(\mu) & = \frac{\partial}{\partial \tau}[C_t V][jY_{ft}]\mu\\
	& = [C_t V][\mu][-jY_{ft}][\tau]^{-1}\\
	& = I_{f_{\tau\theta}}(\mu) = I_{f_{\tau\theta}}(\mu)^T\\
I_{f_{\theta\theta}}(\mu) & = \frac{\partial}{\partial \theta}[C_t V][jY_{ft}]\mu\\
	& = [C_t V][\mu][-Y_{ft}]
\end{align*}

Similar derivations hold for $I_t$:

\begin{align*}
I_{t_{\Theta\tau}}(\mu) & = \frac{\partial}{\partial \tau}(I_{t_\Theta}^T \mu)\\
	& = \frac{\partial}{\partial \tau}\Big(j [V](C_f^T [Y_{tf}] + C_t^T [Y_{tt}]) \mu\Big)\\
	& = \frac{\partial}{\partial \tau}\Big(j [V](C_f^T [\mu] Y_{tf} + C_t^T [\mu]Y_{tt})\Big)\\
	& = j [V]C_f^T [\mu] [-Y_{tf}][\tau]^{-1}\\
I_{t_{\Theta\theta}}(\mu) & = \frac{\partial}{\partial \theta}\Big(j [V](C_f^T [\mu] Y_{tf} + C_t^T [\mu]Y_{tt})\Big)\\
	& = j [V]C_f^T [\mu][-jY_{tf}]\\
I_{t_{\V\tau}}(\mu) & = \frac{\partial}{\partial \tau}(I_{t_\V}^T \mu)\\
	& = \frac{\partial}{\partial \tau}\Big([\V]^{-1}[V](C_f^T [Y_{tf}] + C_t^T [Y_{tt}]) \mu\Big)\\
	& = \frac{\partial}{\partial \tau}\Big([\V]^{-1}[V](C_f^T [\mu] Y_{tf} + C_t^T [\mu]Y_{tt})\Big)\\
	& = [\V]^{-1}[V](C_f^T [\mu] [-Y_{tf}][\tau]^{-1}\\
I_{t_{\V\theta}}(\mu) & = \frac{\partial}{\partial \theta}\Big([\V]^{-1}[V](C_f^T [\mu] Y_{tf} + C_t^T [\mu]Y_{tt})\Big)\\
	& = [\V]^{-1}[V]C_f^T [\mu][-jY_{tf}]
\end{align*}
\begin{align*}
I_{t_{\tau\Theta}}(\mu) & = I_{t_{\Theta\tau}}(\mu)^T\\
I_{t_{\tau\V}}(\mu) & = I_{t_{\V\tau}}(\mu)^T\\
I_{t_{\theta\Theta}}(\mu) & = I_{t_{\Theta\theta}}(\mu)^T\\
I_{t_{\theta\V}}(\mu) & = I_{t_{\V\theta}}(\mu)^T\\
\\
I_{t_{\tau\tau}}(\mu) & = \frac{\partial}{\partial \tau}([\tau]^{-1}[-Y_{tf}][C_fV]\mu)\\
	& = \frac{\partial}{\partial \tau}(-[\tau]^{-1}[C_fV][\mu]Y_{tf})\\
	& = 2[\tau]^{-2}[C_fV][\mu][Y_{tf}])\\
I_{t_{\tau\theta}}(\mu) & = \frac{\partial}{\partial \theta}(-[\tau]^{-1}[C_fV][\mu]Y_{tf})\\
	& = - [\tau]^{-1}[C_fV][\mu][-jY_{tf}]\\
I_{t_{\theta\tau}}(\mu) & = \frac{\partial}{\partial \tau}[C_f V][-jY_{tf}]\mu\\
	& = [C_f V][\mu][jY_{tf}][\tau]^{-1}\\
	& = I_{t_{\tau\theta}}(\mu)^T\\
I_{f_{\theta\theta}}(\mu) & = \frac{\partial}{\partial \theta}[C_f V][-jY_{tf}]\mu\\
	& = [C_f V][\mu][-Y_{tf}]
\end{align*}

\bibliographystyle{plain}
\bibliography{biblio}

\end{document}